\documentstyle[12pt,amscd,amssymb,righttag]{amsart}
\let\Pp=\P

\newtheorem{thm}{Theorem}[section]

\newtheorem{lemma}[thm]{Lemma}
\newtheorem{prop}[thm]{Proposition}

\theoremstyle{definition}
\newtheorem{dfn}[thm]{Definition}

\theoremstyle{remark}

\numberwithin{equation}{section}

\newcommand{\C}{{\Bbb C}}

\newcommand{\OO}{{\scrpt O}}
\let\cross=\times
\let\tensor=\otimes

\newcommand{\Hilb}{\operatorname{Hilb}}
\newcommand{\Ext}{\operatorname{Ext}}

\newcommand{\ch}{\operatorname{ch}}
\newcommand{\rk}{\operatorname{rk}}
\newcommand{\I}{{\cal I}}
\newcommand{\M}{{\cal M}}

\newcommand{\compo}{\raise2pt\hbox{$\scriptscriptstyle\circ$}}
\renewcommand{\leq}{\leqslant}
\renewcommand{\geq}{\geqslant}
\newcommand{\lra}{\longrightarrow}
\newcommand{\lRa}[1]{\>{\buildrel {#1}\over\longrightarrow}\>}

\def\rep#1{\bysame}
\def\bib[#1]#2<#3>#4|#5(#6){\bibitem{#1}
{\sc #2},\ #3,\ {\it #4},\ {\bf #5}\ (#6)}
\def\tbib[#1]#2<#3>#4.{\bibitem{#1} {\sc #2},\ {\it #3},\ #4.}
\def\bibt[#1]#2<#3>#4|#5(#6){\bibitem{#1}
{\sc #2},\ #3,\ in\ {\it #4}\ ed.\ {#5},\ (#6)}
\def\bibit[#1]#2<#3>{\bibitem{#1} {\sc#2},\ #3}
\def\prebib[#1]#2<#3>{\bibitem{#1} {\sc #2},\ {\it #3},\ Preprint}
\def\toabib[#1]#2<#3>{\bibitem{#1} {\sc #2},\ #3,\ to appear}

\let\scrpt=\cal
\let\o=\operatorname
\def\R{{\bold R}}

\def\geq{\geqslant}
\def\leq{\leqslant}
\def\whatx{{\widehat X}}
\def\tf{\widehat{\cal I}_W}
\def\ch{\operatorname{ch}}
\begin{document}

\title[]{HILBERT SCHEMES OF POINTS ON SOME \\
K3 SURFACES AND GIESEKER STABLE BUNDLES}
\author[]{Ugo Bruzzo$\Pp$}
\author[]{Antony Maciocia\S}
\address{$\Pp$\thinspace International School for Advanced Studies\\
Via Beirut 2-4, 34014 Miramare, Trieste, Italy.}
\email{bruzzo@@sissa.it}
\address{\S\thinspace Department of Mathematics and Statistics\\
The University of Edinburgh\\
The King's Buildings\\ Mayfield Road\\ Edinburgh, EH9 3JZ, UK.}
\email{ama@@maths.ed.ac.uk}
\thanks{Research partly supported by the National Group for Mathematical
Physics of the Italian Research Council, by the Italian Ministry for
University and Research through the research project
`Geometria  delle  variet\`a
differenziabili,' the Engineering and Physical Sciences
Research Council of United Kingdom and the Edinburgh Mathematical Society.
The authors would like to thank the Vector Bundle's group of Europroj.
The second author is a member of the VBAC group of Europroj.}
\date{\today}

\maketitle

\begin{abstract}
By using a Fourier-Mukai transform for sheaves on K3 surfaces we
show that for a wide class of K3 surfaces
$X$ the Hilbert schemes $\Hilb^n(X)$ can be identified for all $n\geq
1$ with moduli spaces of Gieseker stable vector bundles on $X$. We also
introduce a new Fourier-Mukai type transform for such surfaces.
\end{abstract}

\section*{Introduction}
Let $(X,H)$ be a polarized K3 surface over $\C$ which also carries a
divisor $\ell$ such that
\begin{equation}
H^2=2\,,\qquad\ell^2=-12\,,\qquad H\cdot\ell=0.\label{req}\end{equation}
One must also include a technical condition which can be expressed in the form
$H^0(\OO(\ell+2H))=0$. This will hold generically.
There is an 18-dimensional family of such K3 surfaces which are called
`reflexive' in \cite{BBHM}; these include generic Kummer surfaces.
For any sheaf $\cal E$ on $X$ we denote its Mukai vector by $v(\cal E)$.
Recall that $v({\cal E})=(\rk\cal E,c_1(\cal E),s(\cal E))$, where
$$s(\cal E)=\rk\cal E+\o{ch}_2(\cal E)\,.$$
For such surfaces it can be shown that the moduli space
$\whatx=\cal M(2,\ell,-3)$ of Gieseker stable sheaves $\cal E$ on $X$ with
Mukai vector
$v(\cal E)=(2,\ell,-3)$ is isomorphic to $X$, and  is formed by locally-free
$\mu$-stable sheaves. The natural polarization on $\hat X$ is $\hat H$ which
one can show is given by $2\ell+5H$ if we identify $\hat X$ with $X$. There is
a divisor $\hat\ell$ on $\hat X$ which plays the role of $\ell$ and one can
show that it is given by $-5\ell-12H$ if we identify $\hat X$ with $X$.

In this paper we show the following result.
\begin{thm}
For any $n\geq 1$, the Hilbert scheme
$\Hilb^n(X)$ of length-$n$ zero-di\-men\-sion\-al subschemes of $X$ is
isomorphic, as an algebraic variety, to the moduli space ${\cal
M}_n={\cal M}(1+2n,-n\hat\ell,1-3n)$ of Gieseker stable bundles on
$X$ with respect to $\hat H$.\label{th:main}
\end{thm}
Notice that the theorem implies that all Gieseker stable sheaves are
locally-free. It will also follow that they are never
$\mu$-stable with respect to $\hat H$. In particular, the underlying smooth
bundles do not carry anti-self-dual connections with respect to the K\"ahler
metric associated to $\hat H$.

The isomorphism $\Hilb^n(X)\simeq\cal M_n$ will be
established by using a Fourier-Mukai (FM) transform of K3 surfaces (see
\cite{BBH}). This FM transform preserves
the natural complex symplectic structures of the moduli spaces
\cite{Mac} and so $\Hilb^n(X)$ and $\cal M_n$ are isomorphic as
complex hyperk\"ahler manifolds as well.

Our results should be compared with those given in \cite{Z}, where for a
general polarized K3 surface $(X,H)$  a birational map
$\cal M(2,0,-1-n^2H^2)\to \Hilb^{2n^2H^2+3}(X)$ is constructed. Further
birational identifications can be found in \cite{BC}.
A similar version of Theorem 1 can be
found in \cite{Mac2} for the case of an abelian surface. Our theorem is
consistent with the results of \cite{GH} which show that moduli spaces
of semi-stable torsion-free sheaves have the same Hodge numbers as Hilbert
schemes of points of the K3 surface.

\section{Fourier-Mukai Transforms for K3 Surfaces}
Let us recall some definitions and properties related to the FM transform on
$X$ introduced in \cite{Muk3}.
Let $\cal Q$ denote the universal sheaf on $X\times\whatx$ normalized by
$\R\hat\pi_*({\cal Q})\cong\OO_{\whatx}[-1]$. This gives rise to a
Fourier-Mukai
transform
$$\R\Phi({\cal E})=\R\hat\pi_*(\pi^*{\cal E}\tensor{\cal Q}).$$
In \cite{BBH} this is shown to give an equivalence of derived categories
between the derived category $D(X)$ of complexes of coherent sheaves on $X$ and
$D(\whatx)$. The inverse is given by
$$\R\hat\Phi({\cal E})=\R\pi_*(\hat\pi^*{\cal E}\tensor{\cal Q}^*)$$
up to a shift of complexes.

\begin{dfn}
We say that a sheaf $\cal E$ on $X$ is
IT$_k$ if $H^j(X,\cal E\otimes\cal Q_\xi)=0$ for $j\neq k$, where
$\cal Q_\xi=\cal Q\vert_{X\times\{\xi\}}$ for $\xi\in\whatx$. We say that
$\cal E$ is WIT$_k$ if $R^j\hat\pi_\ast(\pi^\ast\cal
E\otimes\cal Q)=0$ for $j\neq k$,
where $\pi$ and $\hat\pi$ are the projections of
$X\times\whatx$ onto the two factors. Obviously, any IT$_k$ sheaf is
WIT$_k$.
\end{dfn}

For any WIT$_k$ sheaf $\cal E$ on $X$ we denote
its FM transform as the sheaf on
$\whatx$
$$\widehat{\cal E}=R^k\Phi({\cal E}).$$
Note that $\cal E$ is IT$_k$ if and only if it is WIT$_k$ and
its FM transform is locally free.

The main properties of this FM transform are summarized as follows.
\begin{prop}
Let $\cal E$ be a sheaf on $X$. Then the Chern character of the transform of
$\cal E$ is given by
\begin{align*}
\ch_0&=-\ch_0({\cal E})+c_1({\cal E})\cdot\ell+2\ch_2({\cal E})\\
\ch_1&=-c_1-(c_1\cdot H+5\ch_2({\cal E}))\hat\ell
+(c_1\cdot\ell-2c_1\cdot H)\hat H\\
\ch_2&=-5\ch_2({\cal E})-2c_1\cdot\ell,
\end{align*}
where $c_1=c_1({\cal E})$.
{}From this it follows that $\deg\R\Phi{\cal E}=-\deg\cal E$ and
$\chi(\R\Phi{\cal E})=-\chi({\cal  F})$.
\label{topinv}\end{prop}

\begin{prop} {\rm (Invertibility)} Let $\cal E$ be a WIT$_k$ sheaf on   $X$.
Then its FM transform $\widehat{\cal E}=R^k\hat\pi_\ast(\pi^\ast\cal
E\otimes\cal Q)$ is a WIT$_{2-k}$ sheaf on $\whatx$, whose
inverse FM transform $R^{2-k}\pi_\ast(\hat\pi^\ast\widehat{\cal
E}\otimes{\cal Q}^\ast)$ is isomorphic to $\cal  F$.
\end{prop}

\section{The isomorphism between the Hilbert scheme and the moduli space}
We give a definition that extends to K3 surfaces the notion of
homogeneous bundle on an abelian variety (cf.\ \cite{Muk1}).
\begin{dfn}
A coherent sheaf $\cal F$ on $\whatx$ is {\em quasi-homogeneous} if it
has a filtration by sheaves of the type $\widehat{\OO_Y}$,
where the $Y$'s are zero-dimensional subschemes of $X$, so that
the associated grading is of the form $\oplus_k\cal Q_{p_k}$,
with $p_k\in X$.
\end{dfn}
Let $W$ be a zero-dimensional subscheme of $X$;
the structure sheaf $\OO_W$ is IT$_0$,
and its FM transform $\widehat\OO_W$ is
a quasi-homogeneous locally free sheaf.
$\widehat{\OO}_W$ is $\mu$-semistable, due to the following result.
\begin{prop}
Let $p\in X$. Any nontrivial extension
\begin{equation}
0 \lra \cal Q_p \lra \cal F \lra \cal Q_p \lra 0 \label{ext}
\end{equation}
is $\mu$-semistable. Any destabilizing $\mu$-semistable
subsheaf of $\cal F$ is isomorphic to $\cal Q_p$.
\label{prop1}\end{prop}
\begin{pf}
The first part is standard. For the second part one observes that any
torsion-free destabilising sheaf ${\cal Q}'$ of ${\cal F}$ which is
$\mu$-semistable must have Chern character $(2,\ell,-5)$ because both
${\cal Q}'$ and ${\cal F}/{\cal Q}'$ must satisfy the Bogomolov inequality.
Then it must also be locally-free by the
Bogomolov inequality applied to ${\cal Q}'{}^{**}$.
\end{pf}
\begin{lemma}
The ideal sheaf $\cal I_W$ is IT$_1$.
\end{lemma}
\begin{pf}
Let $\xi\in\whatx$. Then
$H^0(X,\cal I_W\otimes\cal Q_\xi)\hookrightarrow H^0(X,\cal
Q_\xi)=0$ because $\cal Q_\xi$ is $\mu$-stable, and
$$H^2(X,\cal I_W\otimes\cal Q_\xi)^\ast\simeq
\o{Ext}^0(\cal I_W\otimes\cal Q_\xi,\OO_X)\simeq
\o{Hom}(\cal I_W,\cal Q_\xi^\ast)=0$$
since $\cal Q_\xi$ is locally free.
\end{pf}
So by applying the FM transform to the sequence
$$ 0 \lra \cal I_W \lra \OO_X \lra \OO_W \lra 0$$
we get
\begin{equation}
0 \lra \widehat\OO_W \lra \widehat{\cal I}_W \lra \OO_{\widehat
X} \lra 0.\label{due}
\end{equation}

\begin{prop}
The FM transform $\widehat{\cal I}_W$ is Gieseker stable.
\end{prop}
\begin{pf}
Since $\o{ch}(\tf)=(1,0,-n)$, by the formulas in Proposition \ref{topinv} we
obtain
$$\rk\tf=1+2n,\quad\ch_2\tf=-5n,\quad\tilde p(\tf)=\frac{\chi(\tf)}{\rk\tf}=
\frac{2-n}{1+2n}>-\frac12.$$
Let $\cal A$ be a destabilizing subsheaf of $\tf$, that we may assume to be
Gieseker stable with a torsion-free quotient. Then we have $\tilde p(\cal
A)\geq\tilde p(\tf)>-\frac12$. Let $f$ denote the composite ${\cal
A}\to\tf\to\OO_{\whatx}$.

There are two cases:\\
Case (i) $f=0$. Then there is a map $\cal A\to\widehat\OO_W$. Let
$g_k\colon\cal A\to\cal Q_{p_k}$ be the composition of this map with the
canonical projection onto $\cal Q_{p_k}$.  Since $\tilde p(\cal
A)>-\frac12=\tilde p(\cal Q_{p_k})$ and both sheaves are Gieseker stable, we
obtain $g_k=0$ for all $k$, which is absurd.\\
Case (ii) $f\ne 0$. We divide this into two further cases: $\rk\cal A=1$ and
$\rk\cal A>1$.

If $\rk\cal A=1$ we have $\cal A^\ast\simeq\OO_{\whatx}$;
hence the sequence (\ref{due}) splits, which contradicts the inversion
theorem $\widehat{\tf}\simeq\cal I_W$.

If $\rk\cal A>1$ we consider the exact sequences
$$ 0 \lra \cal K_1 \lra \cal A \lRa{h} \cal B \lra 0\quad\mbox{and}\quad
0 \lra \cal B \lra \OO_{\whatx} \lra \cal K_2 \lra 0,$$
where $\rk{\cal K}_2=0,1$. If $\rk{\cal K}_2=1$ then $\cal B=0$,
i.e.\ $f=0$ which is absurd, so that $\rk{\cal K}_2=0$, and $\cal B$ has rank
one. We have an exact commuting diagram
\begin{equation}
\begin{CD}
@. 0 @. 0 @. 0 @.\\
@. @AAA @AAA @AAA @.\\
0 @>>> {\cal K}_3 @>>> {\cal K}_4 @>>> {\cal K}_2 @>>> 0 \\
@. @AAA @AAA @AAA @.\\
0 @>>> \widehat\OO_W @>>> \tf @>>> \OO_{\whatx} @>>> 0 \\
@. @AgAA @AAA @AAh'A @.\\
0 @>>> {\cal K}_1 @>>> {\cal A} @>h>> \cal B @>>> 0\\
@. @AAA @AAA @AAA @.\\
@. 0 @. 0 @. 0 @.
\end{CD}\label{diagtwo}
\end{equation}
with $\mu({\cal K}_1)=0$, $0<\rk{\cal K}_1<2n$ and $f=h'\compo h$.

If $n=1$ then $\widehat\OO_W$ is $\mu$-stable, but this is a contradiction.
For $n>1$, we may assume that $\cal K_1$
is $\mu$-semistable so that it is a direct summand of $\widehat\OO_W$.
Then $\cal K_3$ is locally free and $\rk\cal K_1\geq 2$. Moreover,
$\mu(\cal B)\leq 0$ because
$\cal B$ injects into  $\OO_{\whatx}$, and $\mu(\cal B)\geq 0$ because
$\mu(\cal K_1)\leq 0$. Then $\mu(\cal B)= \mu(\cal K_1)=0$.
Since $\cal K_3$ is locally free the support of $\cal K_2$ is not
zero-dimensional. So $\mu(\cal B)= 0$ implies $\cal K_2=0$ and $\cal
K_3\simeq\cal K_4$.

Finally, we consider the middle column in (\ref{diagtwo}). The sheaf $\cal A$
has rank greater than 2, and is Gieseker stable, so that it is IT$_1$.
But $\tf$ is WIT$_1$ while $\cal K_4$ is WIT$_2$. Then ${\cal
A}\simeq\tf$, but this is a contradiction.
\end{pf}
Note that $\tf$ is never $\mu$-stable because (\ref{due}) destabilizes it.

Let $\cal M_n$ be the moduli space $\cal M(1+2n,-n\hat\ell,1-3n)$ of
Gieseker stable sheaves on $\whatx$. The previous construction yields a map
$\Hilb^n(X)\to\cal M_n$. This map is algebraic because the
Fourier-Mukai transform is functorial and so preserves the Zariski tangent
spaces (see \cite{Mac}). Another way to see this is to observe that the
Fourier-Mukai transforms give a natural isomorphism of moduli functors and so
give rise to an isomorphism of (coarse or fine) moduli schemes.\footnote{We
would like to thank Daniel Hern\'andez Ruip\'erez for this observation.}
We shall now show that the Fourier-Mukai transform is a surjection up to
isomorphism.

\begin{lemma}
Any element $\cal F\in\cal M_n$ is WIT$_1$.
\end{lemma}
\begin{pf} Since $\tilde p(\cal F)>-\frac12$ and $\tilde p(\cal Q_p)=-\frac12$
there
is no map $\cal F\to\cal Q_p$. This means that $H^2(\whatx,\cal
F\otimes\cal Q^\ast_p)=0$.

We consider now nonzero morphisms $\cal Q_p\to\cal F$. Any such
map is injective; otherwise it would factorize through a
rank-one torsion-free sheaf $\cal B$ with $\mu(\cal B)>0$
(because $\cal Q_p$ is $\mu$-stable) and $\mu(\cal B)\leq 0$
(because $\cal F$ is $\mu$-semistable), which is impossible.
Then $\cal Q_p$ is a locally free element of a Jordan-Holder filtration of
$\cal F$. Since any such filtration has only a finite number of
terms, and the associated grading $\o{gr}(\cal F)^{\ast\ast}$
is unique, there is only a finite number of $p$'s giving rise to
nontrivial morphisms, i.e.\ $\o{Hom}(\cal Q_p,\cal F)\simeq
H^0(X,\cal F\otimes\cal Q_p^\ast)$ does not vanish only for a
finite set of points $p$. This suffices to prove that $\cal F$ is WIT$_1$ due
to Proposition 2.26 of \cite{Muk3}.
\end{pf}

\begin{prop}
The FM transform $\widehat{\cal F}$ of $\cal F$ is torsion-free.
\end{prop}
\begin{pf}
Let $\cal T$ be the torsion subsheaf of $\widehat{\cal F}$, so
that one has an exact sequence
\begin{equation}
0 \lra \cal T \lra \widehat{\cal F} \lra \cal G \lra 0.\label{tor}
\end{equation}
Since $\cal T$ is supported at most by a divisor, and
$\widehat{\cal F}$ is WIT$_1$, the sheaf $\cal T$ is WIT$_1$ as
well.
Moreover $\deg(\cal T)\geq 0$. If $\deg\cal T=0$ then $\cal T$ is IT$_0$,
i.e.\ $\cal T=0$.

Hence, we assume $\deg(\cal T)> 0$. The rank-one sheaf $\cal G$ is
torsion-free and, by imbedding it into its double dual, we see that
it is IT$_1$. Then, applying $\R\hat\Phi$ to
(\ref{tor}), we get
$$ 0 \lra \widehat{\cal T} \lra \cal F \lra \widehat{\cal G}
\lra 0.$$
Since $\cal F$ is $\mu$-semistable we see that
$$\deg\cal T=\deg  \widehat{\cal T} \leq 0,$$
which is a contradiction.
\end{pf}
Now the Chern character of $\widehat{\cal F}$ is $(1,0,-n)$, so
that it is the ideal sheaf of a zero-dimensional subscheme of
$X$ of length $n$. We have therefore shown that the Fourier-Mukai transform
surjects as a map $\Hilb^nX\to\M_n$.
The inversion theorem for $\R\Phi$ therefore implies that the transform gives
an isomorphism of smooth varieties. This establishes Theorem \ref{th:main}.

\section{Another Fourier-Mukai Transform}

We shall now show that $\M_1$ gives rise to another FM transform.
\begin{lemma} $\dim H^1({\cal Q})=1$.
\end{lemma}
\begin{pf}
This follows immediately from the degeneration of the Leray spectral sequence
applied to $\hat\pi$ and the fact that $\R\hat\pi_*({\cal
Q})=\OO_{\whatx}[-1]$.
\end{pf}
This lemma shows that there is a unique extension
\begin{equation}
0\lra{\cal Q}\lra{\cal E}\lra\OO_{X\cross\whatx}\lra0.\label{psiE}
\end{equation}
Both of the restrictions of ${\cal E}$ to the factors of $X\cross\whatx$ give
families of Gieseker stable bundles.
Then the general theory of FM transforms (see \cite{Mac2}) implies that ${\cal
E}$ gives rise to an FM transform which we denote by $\R\Psi$.
\begin{prop}The FM transform $\R\Psi$ and its inverse $\R\hat\Psi$ satisfy the
following:
\begin{enumerate}
\item $\R\Psi\OO_X\cong\OO_{\whatx}[-2]$,
\item $\ch(\R\Psi{\cal F})=\chi({\cal F})\ch(\OO_{\whatx})
+\ch(\R\Phi{\cal F})$, where ${\cal F}$ is any sheaf on $X$ and
\item $\R\hat\Psi\I_p\cong{\cal Q}^*_p[-1]$, for $p\in\whatx$.
\end{enumerate}
\end{prop}
\begin{pf}
Apply $\R\hat\pi_*$ to the short exact sequence \ref{psiE} to obtain a long
exact sequence. Note that $\R\hat\pi_*\OO_{X\cross\whatx}$ is concentrated in
the 0th and 2nd positions and $\R\Phi{\cal Q}=\OO_{\whatx}[-1]$. Then part (1)
follows immediately.

For the second part just twist \ref{psiE} by $\pi^*{\cal F}$ and apply
$\R\hat\pi_*$. Then the formula follows from the facts that the Chern
character is additive with respect to triangles in $D(\whatx)$ and
$$\R\hat\pi_*(\pi^*{\cal E})=\R\Gamma(X,{\cal E})\tensor\OO_{\whatx}$$
by the projection formula, where $\Gamma$ denotes the sections functor.

By (1) we have $\R\hat\Psi\OO_{\whatx}=\OO_X$. We also have
$\R\hat\Psi\OO_p\cong{\cal E}^*_p$ from the definition of $\R\hat\Psi$. Then
when we apply $\R\hat\Psi$ to the structure sequence of $p$ we obtain the
short exact sequence
$$0\lra \OO_X\lra{\cal E}^*_p\lra R^1\hat\Psi\I_p\lra0.$$
This is just the dual of (\ref{due}). This completes the proof.
\end{pf}
Note that $\I_W$ does not satisfy WIT with respect to $\R\Psi$ but $\I_{\hat
W}$ does satisfy WIT$_1$ with respect to $\R\hat\Psi$ and the transform is
just the dual of the corresponding quasi-homogeneous bundle.

\section{Concluding remarks}
Since both the moduli spaces and the punctual Hilbert schemes have complex
symplectic structures which are given by the cup product on $\Ext^1(E,E)$ the
FM transforms will preserve the symplectic structures and so are complex
symplectic isomorphisms. It follows immediately that the spaces are
hyperk\"ahler isometric as well.

Theorem \ref{th:main} has several immediate consequences which we can state in
the
following theorem.
\begin{thm} Let $X$ be a reflexive K3 surface.
\begin{enumerate}
\item The moduli space $\M_n$ of Gieseker stable sheaves on $X$ is connected
and projective. All points of $\M_n$ are locally-free.
\item $\M_n$ contains no $\mu$-stable sheaves and so the moduli space of
irreducible U$(2n+1)$-instantons, with fixed determinant $\OO(\hat\ell)^{-n}$
and second Chern character $-5n$, is empty.
\item The moduli space of all instantons with this type is isomorphic to the
$n^{th}$ symmetric product $S^nX$.
\end{enumerate}
\end{thm}
For the last part one uses the fact that any $\mu$-semistable sheaf of the
given Chern character admits a surjection to $\OO_X$ and so fits into a
sequence of the form (\ref{due}).

\end{document}